\documentclass[a4paper,12pt]{article}
\usepackage[T2A]{fontenc}                      
\usepackage[cp1251]{inputenc}           
\usepackage[all]{xy}
\usepackage{amssymb}
\usepackage{cite}
\usepackage{cmap}
\usepackage{latexsym}
\usepackage{enumerate}
\usepackage{amsmath, amsthm, amscd, amsfonts, amssymb, graphicx, color}
\usepackage[left=2.5cm,right=2.5cm,top=2cm,bottom=2cm]{geometry}
\usepackage{indentfirst}
\usepackage{array}
\usepackage{bm}
\usepackage{float}

\usepackage{wrapfig}

\usepackage{authblk}

\title{A short note on $K_{n}$-irregular graphs}

\author{Tatiana~Dovzhenok\thanks{Corresponding author. E-mail: \texttt{t.dovzhenok@mail.ru}}}

\affil{\small Research Laboratory ``Mathematics of Hybrid Intelligence Systems'', \\	
Francisk Skorina Gomel State University, Gomel, 246028, Belarus}

\date{} 

\begin{document}	
\maketitle

\begin{abstract}
This note addresses Theorem 12 in the seminal paper ``$F$-degrees in graphs'' by G.~Chartrand et al. (1987), which asserts the existence of $K_n$-irregular graphs for each $n \ge 3$. Although this result is widely cited in the literature as an established fact, the original text lacks a complete proof. We show that for $n = 5$ and all $n \ge 7$, the constructions as presented are not $K_n$-irregular, presumably due to typographical errors. Consequently, referencing this foundational work as containing a rigorous proof for the general case is not entirely accurate.

\end{abstract}

\section{Introduction}
In 1987, G.~Chartrand, K.~S.~Holbert, O.~R.~Oellermann, and H.~C.~Swart~\cite{1} introduced the concept of $F$-degrees in graphs. For a fixed graph $F$, the \emph{$F$-degree} of a vertex $v$ in a graph $G$ is defined as the number of subgraphs of $G$ isomorphic to $F$ that contain~$v$. A graph $G$ is called \emph{$F$-irregular} if all its vertices have pairwise distinct $F$-degrees.

Let $K_n$ denote the complete graph of order $n$. In the case of such graphs, this existence statement is formulated in the original paper as Theorem 12:
\begin{quote}
	\textbf{Theorem 12.} \emph{For every integer $n \ge 3$, there exists a nontrivial $K_n$-irregular graph.}
\end{quote}

This assertion has been traditionally cited as proven in~\cite{1} (see, e.g., \cite{2,3,4}), a practice also followed in the present author's papers, \emph{inter alia}~\cite{5,6}. Since the original work remains publicly unavailable in digital form, a direct verification of its results is difficult.

In Section 2 of this note, we recall the classical constructions proposed for Theorem~12 in~\cite{1} and demonstrate that they do not possess the required $K_n$-irregularity for $n = 5$ and $n \ge 7$. Combined with the lack of an explicit irregularity analysis in the original text, this leaves the theorem unproven for these values of $n$. Completing the general case requires either correcting these constructions or designing new ones.

The identified gap becomes particularly relevant in light of a recent preprint by J.~A.~Schreib~\cite{7}, which claims a proof of the existence conjecture on $F$-irregular graphs (Chartrand et~al.~\cite{1}). Therein, the case of complete graphs is taken as fully resolved, relying solely on~\cite{1}.

Finally, infinite families of $K_n$-irregular graphs for each $n \ge 3$ were recently reported by the author and A.~Filuta~\cite{8} as part of the proof of the strong conjecture about \mbox{$F$-irregular} graphs (Dovzhenok et al.~\cite{5}).

\section{Structural Analysis of Constructions for Theorem 12}

To justify Theorem 12 for $n \ge 5$ and $n \neq 6$, the authors in~\cite{1} propose, without $K_n$-degree verification, two block constructions based on the parity of $n$.

For odd $n \ge 5$, the graph $G$ is defined by
\begin{align*}
	V(G) &= \{a_1, a_2, a_3, b_1, b_2, \dots, b_{n+1}, c_1, c_2, \dots, c_{n-1}, d\}, \\
	E(G) &= \{a_i b_j \mid i \le j, \, 1 \le i \le 3, \, 1 \le j \le n+1\} \\
	&\quad \cup \{b_j b_l \mid 1 \le j < l \le n+1\} \\
	&\quad \cup \{b_j c_k \mid 1 \le k \le n-1, \, 1 \le j \le n, \, k \le j+1\} \\
	&\quad \cup \{c_k c_m \mid 1 \le k < m \le n-1\} \\
	&\quad \cup \{b_j d \mid 2 \le j \le n-1\} \cup \{a_3 c_{n-1}, a_2 d\}.
\end{align*}

For even $n \ge 8$, the graph $G$ is given by
\begin{align*}
	V(G) &= \{a_1, a_2, a_3, b_1, b_2, \dots, b_{n+1}, c_1, c_2, \dots, c_{n-1}, d\}, \\
	E(G) &= \{a_i b_j \mid i \le j, \, 1 \le i \le 3, \, 1 \le j \le n+1\} \\
	&\quad \cup \{b_j b_l \mid 1 \le j < l \le n+1\} \\
	&\quad \cup \{b_j c_k \mid 1 \le k \le n-1, \, 1 \le j \le n, \, 1 \le k \le j+1\} \\
	&\quad \cup \{c_k c_m \mid 1 \le k < m \le n-1\} \\
	&\quad \cup \{b_j d \mid 2 \le j \le n-3\} \cup \{a_3 c_{n-1}, a_1 d, a_2 d, b_{n+1} d\}.
\end{align*}

In each such graph $G$, the vertices $c_1$ and $c_2$ have identical closed neighborhoods:
\begin{equation*}
	N[c_1] = N[c_2] = \{b_1, b_2, \dots, b_{n}\} \cup \{c_1, c_2, \dots, c_{n-1}\}.
\end{equation*}

Due to this symmetry, the $K_n$-degrees of the vertices $c_1$ and $c_2$ in $G$ are equal, which contradicts the required $K_n$-irregularity, rendering these specific constructions invalid for the proof of Theorem 12.

\end{document}